\documentclass{amsart}

\usepackage{etoolbox} 

\newtheorem{defn}{Definition}

\newtheorem{rem}{Remark}

\usepackage[pdfauthor={Richard J. Mathar},colorlinks=true,bookmarksopen]{hyperref}

\usepackage{orcidlink}
\begin{document}

\title[Roots of Zernike Polynomials]{Third Order Newton's Method for Zernike Polynomial Zeros}

\author{Richard J. Mathar \orcidlink{0000-0001-6017-6540}}
\urladdr{https://mathar.www3.mpia-hd.mpg.de}
\email{mathar@mpia-hd.mpg.de}
\address{Max-Planck Institute for Astronomy, K\"onigstuhl 17, 69117 Heidelberg, Germany}
\thanks{Supported by the NWO VICI grant 639.043.201 to A. Quirrenbach,
``Optical Interferometry: A new Method for Studies of Extrasolar Planets.''}


\date{\today}

\begin{abstract}
The Zernike radial polynomials are a system of orthogonal polynomials
over the
unit interval with weight $x$.
They are used as basis functions in optics to
expand
fields over the cross section of circular pupils.
To calculate the roots of Zernike polynomials,
we optimize the generic iterative numerical
Newton's Method that iterates on zeros of functions with third order convergence.
 The technique
is based on rewriting the polynomials as Gauss Hypergeometric Functions, reduction of
second order derivatives to first order derivatives, and evaluation of some
ratios of derivatives by terminating continued fractions.

A PARI program and a short table of zeros
complete up
to polynomials of 40th order are included.
\end{abstract}

\maketitle
\section{Classical Orthogonal Polynomials: Hofsommer's Method} 

The generic third order Newton's Method---also known as Halley's method---to compute roots $f(x)=0$ numerically
improves solutions $x_i\rightarrow x_{i+1}=x_i+\Delta x$ iteratively, starting
from initial guesses, via computation of corrections
\begin{equation}
\Delta x = -\frac{f(x)}{f'(x)}/\left(1-\frac{f(x)}{2f'(x)}\,\frac{f''(x)}{f'(x)}\right)
\label{eq.New}
\end{equation}
where $f(x)$, $f'(x)$ and $f''(x)$ are the function and its first and second
derivatives at the current best approximation $x_i$ \cite{GerlachSIAM36,HansenNumerMath27,KalantariJCAM80}.
For some
classes of orthogonal polynomials, $f''/f'$ can be derived
from $f/f'$ \cite{HofsommerMTAC12,Shao,WynnMTAC10}, which means the update
can be done to third order at essentially no additional numerical expense.
If we divide the differential equation of the classical orthogonal polynomials,
for example as tabulated in \cite[22.6]{AS}\cite{LewanowiczApplMath29},
\begin{equation}
h_2(x)f''+h_1(x)f'+h_0(x)f=0,
\label{eq.2deq}
\end{equation}
through $f'$, (\ref{eq.New}) turns into
\begin{equation}
\Delta x = -\frac{f(x)}{f'(x)}/
\left[1+\frac{1}{2h_2(x)}\frac{f(x)}{f'(x)}\left(h_0(x)\frac{f(x)}{f'(x)}+h_1(x)\right)\right].
\end{equation}
Structure relations \cite{MarcellanJCAM200} relate the ratio $f/f'$
to ratios at shifted
indices $n$ as tabulated for example in \cite[22.8]{AS},
\begin{eqnarray}
&& g_2(x)f_{n}'(x)=g_1(x)f_n(x)+g_0(x)f_{n-1}(x); \\
\Rightarrow
&&
\frac{f_n(x)}{f_n'(x)}=\frac{g_2(x)}{g_1(x)+g_0(x)\frac{f_{n-1}(x)}{f_n(x)}}.
\label{eq.fnrat}
\end{eqnarray}
The benefit is that the three-term recurrence equations, in the
notation of \cite[22.7]{AS}
\begin{equation}
a_{1,n-1}f_n(x)=(a_{2,n-1}+a_{3,n-1}x)f_{n-1}(x)-a_{4,n-1}f_{n-2}(x) ,
\end{equation}
lead to terminating continued fraction representations for $f/f'$
\begin{equation}
\frac{f_{n-1}(x)}{f_n(x)}=\frac{a_{1,n-1}}{a_{2,n-1}+a_{3,n-1}x-a_{4,n-1}\frac{f_{n-2}(x)}{f_{n-1}(x)}} .
\label{eq.cfrac}
\end{equation}
This is recursively inserted into the denominator  of (\ref{eq.fnrat}) to lower
the index $n$
until $f_0/f_1$ is reached, which avoids problems with cancellation of digits.

This work here implements this strategy
for the Zernike polynomials, $f=R_n^m$,
namely
(i) fast calculation of $f''/f'$ from $f/f'$,
(ii) calculation of $f/f'$ from terminating continued fractions,
both without evaluation of $f$ or its derivatives via direct methods
like Horner schemes.

\section{Zernike Polynomials: Derivatives and Roots} 
\subsection{Definition}
We define Zernike radial polynomials in the unit Ball of dimension $D$
in Noll's nomenclature
\cite{NollJOSA66,PrataAO28,KintnerJModOpt23,TysonOL7,ConfortiOL8,TangoApplPhys13} for $n\ge 0$, $n-m =0 \pmod 2$, $m\le n$ as
\begin{defn} (Zernike Polynomial)
\begin{eqnarray}
R_n^m(x)
=
\sum_{s=0}^{(n-m)/2}
(-1)^s \binom{\frac{n-m}{2}}{s}\binom{\frac{D}{2}+n-s-1}{\frac{n-m}{2}} x^{n-2s}
.
\label{eq.Rpolyx}
\end{eqnarray}
\end{defn}
\begin{rem}
The normalization for $D>2$ might be chosen differently \cite{MatharSAJ179}.
\end{rem}
Following the original notation, we will not put the upper (azimuth) index $m$ in $R_n^m$---which
is not a power---into parentheses.
The normalization integral is
\begin{equation}
\int_0^1 x^{D-1}\, R_n^m(x) R_{n'}^m(x) dx=\frac{1}{2n+D}\delta_{n,n'}.
\label{eq.ortho}
\end{equation}

\begin{rem}
The inversion of (\ref{eq.Rpolyx}) assembles powers $x^i$ by sums
of $R_{n}^m(x)$, ($i\ge m$, $i-m$ even),
\begin{equation}
x^i \equiv 
\sum_{n=m \pmod 2}^i h_{i,n,m} R_n^m(x);\quad
i-m =0,2,4,6,\ldots.
\end{equation}
By projection on the $R_n^m$ basis these moments are
\begin{multline}
h_{i,n,m}
=(2n+D)
\int_0^1 x^{D+i-1} R_n^m(x) dx
\\
=(2n+D)
\sum_{s=0}^{(n-m)/2}
\frac{(-1)^s}
{n-2s+D+i} 
\binom{\frac{n-m}{2}}{s}\binom{\frac{D}{2}+n-s-1}{\frac{n-m}{2}} 
\\
=\frac{n+D+m}{n+D+i}
\binom{n+\frac{D}{2}}{\frac{n-m}{2}}
{}_3F_2\left(
\begin{array}{c}
-(n-m)/2, -(n+D+i)/2 , 1-(D+n+m)/2 \\
1-(n+D+i)/2 , 1-D/2-n \\
\end{array}
\mid 1
\right).
\end{multline}
This terminating Saalsch\"utzian Hypergeometric Function
has a closed form in terms of $\Gamma$-functions \cite{KoornwinderJCAM99}\cite[(2.3.1.3)]{SlaterHyp}:
\begin{equation}
h_{i,n,m}=
(-)^{(n-m)/2}
(n+D/2)
\frac{
(\frac{m-i}{2})_{(n-m)/2}
}{
(\frac{m+i+D}{2})_{1+(n-m)/2}
}
.
\end{equation}
\end{rem}
The round parenthesis with nonnegative lower indices in this equation are the finite products
known as shifted factorials:
\begin{defn}(Pochhammer Symbol)
\begin{equation}
(a)_0\equiv 1;\quad (a)_n=a(a+1)(a+2)\cdots (a+n-1)=\Gamma(a+n)/\Gamma(a);\quad n\ge 0.
\label{eq.Poch}
\end{equation}
\end{defn}

Much
of this work is based on the representation as a terminating
Gaussian Hypergeometric Function,
the product of $x^m$ by a polynomial of degree $n-m$:
\begin{equation}
R_n^m(x)
=
(-1)^{(n-m)/2} \binom{\frac{D+m+n}{2}-1}{\frac{n-m}{2}}
x^m
{}
F\left(
\begin{array}{c}
-(n-m)/2, (D+n+m)/2\\
m+D/2
\end{array}
\mid x^2
\right)
,
\label{eq.RofF}
\end{equation}
with the three ``numerator'' and ``denominator'' parameters
\begin{equation}
a\equiv -(n-m)/2;\quad b\equiv (D+n+m)/2;\quad c\equiv m+D/2,
\end{equation}
and argument
\begin{equation}
z\equiv x^2.
\end{equation}

The initial terms of the power series \eqref{eq.RofF} are
\begin{eqnarray}
R_n^m(x)
=
(-1)^a
{(n+m+D)/2-1\choose (n-m)/2}
x^m
\Big[
1-\frac{(n-m)(D+n+m)}{2(D+2m)}x^2
\nonumber
\\
+\frac{(n-m)(n-m-2)(D+n+m)(D+n+m+2)}{8(D+2m)(D+2m+2)}x^4+\cdots\Big]
.
\label{eq.powxmsplit}
\end{eqnarray}
$R_n^m$ are also Jacobi Polynomials 
\cite[15.4.6,22.5.42,22.5.1]{AS}\cite{ChelyshkovETNA25,TatianJOSA64}:
\begin{multline}
R_n^m(x)
=
(-1)^{(n-m)/2}x^m P_{(n-m)/2}^{(m+1-D/2,0)}(1-2x^2)
\\
=\binom{n+1-D/2}{(n-m)/2} x^m G_{-a}(2+m-D/2,2+m-D/2,x^2)
.
\label{eq.RofP}
\end{multline}
\begin{rem}
Gaussian integration rules for integrals $\int_0^1 x^{D-1}R_n^m(x)f(x) dx$
do not exist because $R_n^m$ changes sign over the integration interval \cite{RothmannMC15}.
(i) \eqref{eq.RofF} suggests to split $R_n^m$ by assigning the factor $x^m$
to the weight such that a Gauss-Legendre integration for moments $x^{D+m-1}$
is engaged and the wiggly remainder of $R_n^m$ multiplied by $f(x)$ is
sampled over the abscissae.
(ii) The representation \eqref{eq.Rpolyx} supports a hybrid term-wise method
that adds the results of $1+(n-m)/2$
Gaussian integrations for moments $x^{D-1+n-2s}$. The disadvantage
here is the need for dense samples of $f(x)$.
(iii) There might be a workaround by developing
rules for
the lifted integrals $\int_0^1 x^{D-1}[1+R_n^m(x)]f(x) dx$ 
following an idea by Denich and Novati \cite{DenichBIT63}.
\end{rem}

\subsection{Recurrence}
Lowering/rising the radial quantum number $n$ by 2 implies increasing/decreasing 
$a$ and $b$ by 1 while keeping $c$ constant.
Contiguous relations (1.27) and (2.1) by Rakha et al.\ 
relate the Gauss Hypergeometric Function to the values with upper parameters 
increased or decreased by 1
\cite[15.2.10]{AS}\cite{RakhaCMA61,ChoEAMJ15}:
$$
(c-a-1)F=(b-a-1)(1-z)F(a^+)+(c-b)F(a^+,b^-);
$$
$$
[(a-b)(a-b-1)(1-z)+b(c-b-1)]F=a(a-b-1)(1-z)F(a^+)+b(c-a)F(a^-,b^+).
$$
Elimination of $F(a^+)$ establishes the format matching $n$-increments and 
$n$-decrements
in \eqref{eq.RofF}:
\begin{multline}
(a-b)[(1+a-b)(1+b-a)(1-z) +a(a-c)+b(b-c)+c-1]F=
\\
-a(c-b)(a-b-1)F(a^+,b^-)
+b(a-b+1)(a-c)F(a^-,b^+).
\end{multline}
Substitution of \eqref{eq.RofF}
$$
F(a,b;c;z) = (-1)^a x^{-m} \frac{1}{\binom{b-1}{-a}}R_n^m
$$
provides the recurrence
\begin{multline}
(1-a)(a-b+1)(a-c)
R_{n+2}^m(x)
=(b-c)(a-b-1)
 (b-1)R_{n-2}^m(x)
\\
-(a-b)[(1+a-b)(1+b-a)(1-z) +a(a-c)+b(b-c)+c-1]R_n^m(x)
;
\end{multline}
\begin{multline}
-(1+\frac{n-m}{2})(1-n-\frac{D}{2})\frac{n+m+D}{2}
R_{n+2}^m(x)
=
\\
\frac{n-m}{2}(1+n+\frac{D}{2})(1-\frac{n+m+D}{2})
R_{n-2}^m(x)
\\
+(n+\frac{D}{2})\left[(1+n+\frac{D}{2})(1-n-\frac{D}{2})(1-x^2) 
+\frac12(n-m)(D+n+m)
+m+\frac{D}{2}-1\right]R_n^m(x)
.
\end{multline}
which generalizes earlier results to $D\ge 2$ \cite{KintnerOA23,PrataAO28}.

\subsection{Derivatives}
Derivatives of (\ref{eq.RofF}) are by the Leibnitz Theorem \cite[3.3.8]{AS}
\begin{eqnarray}
\frac{d}{dx}
R_n^m(x)
&=&
(-1)^{a} { b-1 \choose -a}
\Big[
\frac{d}{dx} x^m F(a,b;c;z)
+x^m
\frac{d}{dx}F(a,b;c;z)\Big];
\label{eq.dRdx}
\\
\frac{d^2}{dx^2}
R_n^m(x)
&=&
(-1)^{a} { b-1 \choose -a}
\Big[
\frac{d^2}{dx^2} x^m F(a,b;c;z)
+
2
\frac{d}{dx} x^m
\frac{d}{dx}
F(a,b;c;z)
\\
&&
\quad\quad
+
x^m
\frac{d^2}{dx^2}F(a,b;c;z)\Big];
\nonumber
\\
\frac{d^3}{dx^3}
R_n^m(x)
&=&
(-1)^{a} { b-1 \choose -a}
\Big[
\frac{d^3}{dx^3} x^m F(a,b;c;z)
+
3
\frac{d^2}{dx^2} x^m
\frac{d}{dx}
F(a,b;c;z)
\label{eq.d3Rdx}
\\
&&
\quad\quad
+
3
\frac{d}{dx} x^m
\frac{d^2}{dx^2}
F(a,b;c;z)
+
x^m
\frac{d^3}{dx^3}F(a,b;c;z)
\Big]
.
\nonumber
\end{eqnarray}
Fa\`a di Bruno's Formula \cite[0.430.2]{GR}\cite{JohnsonAMM109} relegates
the derivatives w.r.t.\ $x$ to derivatives w.r.t.\ $z$,

\begin{eqnarray}
\frac{d}{dx}F(a,b;c;z)
&=&
2xF'(a,b;c;z)
\nonumber
;
\\
\frac{d^2}{dx^2}F(a,b;c;z)
&=&
2F'(a,b;c;z)+4x^2F''(a,b;c;z)
\nonumber
;
\\
\frac{d^3}{dx^3}F(a,b;c;z)
&=&
12x F''(a,b;c;z)+8x^3F'''(a,b;c;z)
\nonumber
 .
\end{eqnarray}
After insertion of these three formulas into (\ref{eq.dRdx})--(\ref{eq.d3Rdx}),
the derivatives of $R_n^m\cong x^mF$ are
\begin{eqnarray}
{R_n^m}'(x)
&\cong&
mx^{m-1}F+2x^{m+1}F'
;
\label{eq.Rpr}
\\
{R_n^m}''(x)
&\cong&
m(m-1)x^{m-2}F+2(2m+1)x^mF'+4x^{m+2}F''
;
\\
{R_n^m}'''(x)
&\cong&
m(m-1)(m-2)x^{m-3}F+6m^2x^{m-1}F'
\\
&&\quad\quad
+12(1+m)x^{m+1}F''+8x^{m+3}F'''
\nonumber
,
\end{eqnarray}
where $\cong$ means the binomial factor and the argument list $(a,b;c;z)$
of the hypergeometric function have not been written down explicitly.
\begin{rem}
These equations write the vector of the derivatives of $R_n^m$ as a lower
triangular matrix with monomials of $x$ multiplied by the vector of
derivatives of $F$. Matrix inversion gives the reciprocal relations:
\begin{eqnarray*}
F &\cong& x^{-m}R ;\\
F' &\cong& -\frac{m}{2}x^{-m-2}R +\frac{1}{2}x^{-m-1}R';\\
F'' &\cong& \frac{m(m+2)}{4}x^{-m-4}R -\frac{1+2m}{4}x^{-m-3}R'+\frac{1}{4}x^{-m-2}R'' ;\\
F''' &\cong& -\frac{m(m+2)(m+4)}{8}x^{-m-6}R +\frac{3(m^2+3m+1)}{8}x^{-m-5}R'
  \\ && -\frac{3(1+m)}{8}x^{-m-4}R''+\frac{1}{8}x^{-m-3}R'''.
\end{eqnarray*}
\end{rem}
Since $R_n^m(x)$ is a polynomial of order $n$, the $(n+1)$st derivatives
equal zero.
Backward elimination of $F$ and its derivatives with the aid of
\cite[15.5.1]{AS}
\begin{equation}
z(1-z)F''(a,b;c;z)+[c-(a+b+1)z]F'(a,b;c;z)=abF(a,b;c;z)
\end{equation}
leads to the analog of (\ref{eq.2deq}),
\begin{multline}
x^2(x^2-1)\frac{d^2}{dx^2}R_n^m(x)
=
\left[nx^2(n+D)-m(D-2+m)\right]R_n^m(x)
\\
+
x\left[D-1-(D+1)x^2\right]
\frac{d}{dx}R_n^m(x).
\label{eq.R2pr}
\end{multline}
This is one special case of differential equations that generate
orthogonal functions \cite{MasjedJMAA325},
and could also be obtained by applying the derivatives of \cite[22.6.1]{AS}
\cite{ElliottMathComp25,DohaJPA35}
to (\ref{eq.RofP}).
The derivative of this reaches out to the third derivatives,
in which $R''$ is reduced to $R$ and $R'$ substituting the previous equation,
\begin{multline}
x^3(x^2-1)^2\frac{d^3}{dx^3}R_n^m(x)
\\
=
\Big\{-n(3+D)(n+D)x^4+\left[(n+m)D^2+(n^2+m^2-n+3m)D-10m+5m^2-n^2\right]x^2
\\
-m(D+1)(D-2+m)\Big\}R_n^m(x)
\\
+
x\Big\{\left[D^2+(3+n)D+n^2+2\right]x^4+\left[-2D^2-(2+n+m)D+6+2m-n^2-m^2\right]x^2
\\
+D^2+D(m-1)-2m+m^2\Big\}\frac{d}{dx}R_n^m(x)
.
\label{eq.R3pr}
\end{multline}

\subsection{Zeros}
\subsubsection{Ratios of Derivatives}

Installation of $f/f'$ in (\ref{eq.New}) progresses by dividing $R_n^m\cong x^mF$
through (\ref{eq.Rpr}),
\begin{equation}
\frac{R_n^m(x)}{{R_n^m}'(x)}
=
\frac{x}{m+2z\frac{F'(a,b;c;z)}{F(a,b;c;z)}}
.
\label{eq.RoverRp}
\end{equation}
The analog of (\ref{eq.fnrat}) is implemented by substituting
\cite[15.2.1]{AS}
\begin{equation}
F'(a,b;c;z)
=\frac{ab}{c}
F(a+1,b+1;c+1;z)
\label{eq.FpoverF}
\end{equation}
in the denominator.
In lieu of (\ref{eq.cfrac}) we find the continued fractions \cite{FrankTAMS81}
\begin{equation}
\frac{F(a,b;c;z)}{F(a+1,b+1;c+1;z)} \equiv
\frac{-bz}{c}+1-
\frac{\frac{(a+1)(c-b)z}{c(c+1)}}{\frac{(a+1-b)z}{c+1}+1-\cdots}
\,
\frac{\frac{(a+2)(c+1-b)z}{(c+1)(c+2)}}{\frac{(a+2-b)z}{c+2}+1-\cdots}
\label{eq.Fcfrac}
\end{equation}
which terminate in our cases since $a$ is a negative integer and $c=a+b$.
This already suffices to implement the standard Newton iteration, i.e., to
approximate (\ref{eq.New}) by $\Delta x = -f(x)/f'(x)$.
Division of (\ref{eq.R2pr}) through ${R_n^m}'(x)$ yields

\begin{equation}
\frac{{R_n^m}''(x)}{{R_n^m}'(x)}
=
\frac{1}{x^2-1}\left[
\left(n(n+D)-\frac{m(D-2+m)}{x^2}\right)
\frac{R_n^{m}(x)}{{R_n^m}'(x)}
+
\frac{D-1-(D+1)x^2}{x}
\right].
\label{eq.R2poverR}
\end{equation}
This is $f''(x)/f'(x)$ of the generic formula, and can be quickly
computed from $R_n^m(x)/{R_n^m}'(x)=f(x)/f'(x)$
of the lower order.

\subsubsection{Initial Guesses}

For $n$ and $m$ fixed, the strategy adopted here is to compute the $(n-m)/2$
distinct roots in $(0,1)$ starting with the smallest, then bootstrapping the others
in naturally increasing order. An approximation to the smallest root is found
by equating the first three terms in the square bracket of (\ref{eq.powxmsplit})
with zero---hoping
that higher powers of $x$ become insignificant for small $x$---and solving the
bi-quadratic equation for $x$. This guess may become unstable for $n$ approximately
larger than 11 in the sense that the Newton iterations converge to another than
this smallest root. Instead, the
simple, heuristic initial guess
\begin{equation}
x\approx \frac{1.46m+2.41}{n+0.46m+1.06}
\label{eq.xinit}
\end{equation}
is used in general, but keeping the solution to the bi-quadratic equation
when this is exact, i.e., in the cases $n-m=2$ or $4$\@.

A shooting method is useful to
produce
an initial estimate for one
root supposed an adjacent one has already been found. The third order Taylor
extrapolation from one root $x$ to the next at $x+\Delta x$ is
\begin{equation}
f(x+\Delta x)\approx
f(x)+\Delta x f'(x)+\frac{(\Delta x)^2}{2!}f''(x)
+\frac{(\Delta x)^3}{3!}f'''(x)\approx 0.
\end{equation}
Division through
$\Delta x f'$
and exploiting
$f(x)=0$
yields a quadratic equation for the approximate distance $\Delta x$ to the next one,
\begin{equation}
1+\frac{\Delta x}{2}\frac{f''(x)}{f'(x)}
+\frac{(\Delta x)^2}{6}\frac{f'''(x)}{f'(x)}\approx 0,
\label{eq.extrap}
\end{equation}
from which the branch $\Delta x>0$ is systematically selected to
start computation of the root adjacent to the previous one.
The two ratios of derivatives are obtained by setting $R_n^m(x)=0$
in (\ref{eq.R2pr}) and (\ref{eq.R3pr}), then dividing both equations through ${R_n^m}'(x)$.
This aim to locate the next root with sufficient accuracy---and to prevent
the Newton's Method to be be
drawn into
the second next root which would
call for more administrative care \cite{EhrlichCACM10}---is the
rationale to look into third derivatives; it might also guide the way
to even higher order Newton's methods employing $f'''$.

\section{Summary} 
The Newton's Method of third order convergence is implemented for
Zernike Polynomials $R_n^m$ by computation of the ratios
${R_n^m}''/{R_n^m}'$ and $R_n^m/{R_n^m}'$ with relay to
the generic formulas of associated, terminating hypergeometric
series. Adding knowledge on the derivative ${R_n^m}'''$, a shooting
method is proposed which generates an initial guess for the adjacent root
from each root found.

\appendix
\section{Table of Roots of Low-order Polynomials}
\label{sec.tabl}
The roots $x_{i,n,m}$ of $R_n^m(x)$ are tabulated 
in the file \texttt{anc/zern2.txt}
for $D=2$,
and in the file \texttt{anc/zern3.txt}
for $D=3$.
Each line contains $n$, then
$m$, then $(n-m)/2=|a|$ values of $x_{i,n,m}$,
the derivative ${R_n^m}'(x_i)$ at the root,
and the weight for the barycentric interpolation 
\cite{BerrutSIAMR46,WangMCom83}.
(The inverse product of the differences to the other roots, ignoring the 
root at $x=0$, i.e., dropping the factor $x^m$ in \eqref{eq.RofF}.)
\begin{rem}
The roots are inert if the normalization factor in \eqref{eq.Rpolyx} were
changed, but the derivatives are not.
\end{rem}
Only the roots $x>0$ are included, and only the
standard parameter range for even, positive values of $n-m$ is considered.

\begin{rem}
For $D=2$, the 
polynomials
$G_{-a}(m+2-D/2,m+2-D/2,y)$ 
in (\ref{eq.RofP}) build
a system of orthogonal polynomials with weight $y^m$ over the unit interval $0\le y\le 1$.
Therefore
the squares $x_{i,n,m}^2=y_{i,n,m}$
are also the abscissae for Gaussian integration of moment $m$
\cite[Tab.\ 25.8]{AS}\cite{FishmanMTAC11,Golub,MatharVixra1303,DuanTAP39}, i.e.,
the $x_{i,n,m}^2$ listed in \texttt{zern2.txt}
 are the abscissae for Gauss-Legendre quadratures for weight $x^m$
on $(n-m)/2$ nodes.
\end{rem}

\section{PARI implementation}
The source code of the PARI interpreter program \texttt{zern.gp}
which computes
the values
tabulated in the \texttt{anc} directory
via \verb+gp -q < zern.gp+ is also listed in the \texttt{anc} directory.
The language
is
similar to C/C++ and has inherent support for
arbitrary precision computation \cite{PARI}.

\texttt{Pochhammer} implements \eqref{eq.Poch}.
\texttt{Hyperg} calculates $F(a,b;c;z)$ for non-positive integer $a$.
\texttt{ZernikeP} implements \eqref{eq.RofF}.
\texttt{ZernikePderiv} calculates ${R_n^m}'(x)$.
\texttt{HypergAugmRatio} implements \eqref{eq.Fcfrac}.
\texttt{HypergRatio} implements \eqref{eq.FpoverF}.
\texttt{ZernikePrratio} implements \eqref{eq.RoverRp}.
\texttt{Zernike2Prratio} implements \eqref{eq.R2poverR}.
\texttt{Zernike3Prratio} implements \eqref{eq.R3pr}.
\texttt{ZernikeRoot} implements \eqref{eq.New}.
\texttt{ZernikeRootEst} implements \eqref{eq.xinit}, but \eqref{eq.powxmsplit}
if $n=m+2$.
\texttt{ZernikeAllRoot} implements a loop with guesses as in \eqref{eq.extrap}.
\texttt{main} loops over $n$ and $m$ to tabulate the zeros up
to a maximum $n$.

\bibliographystyle{amsplain}
\bibliography{all}

\end{document}